\documentclass[11pt]{amsart}
\usepackage{amsmath}
\usepackage{cases}
\usepackage{mathrsfs, euscript}
\usepackage{bbm}
\usepackage{amssymb}
\usepackage{txfonts}
\usepackage{amscd}
\usepackage{amsfonts,latexsym,amsmath,amsthm,amsxtra,mathdots,amssymb,latexsym,mathabx}
\usepackage[all,cmtip]{xy}
\RequirePackage{amsmath} \RequirePackage{amssymb}
\usepackage{color}
\usepackage{colordvi}
\usepackage{multicol}
\usepackage{hyperref}

\def\mod{\mathrm{mod}\,}

\def\sumx{\sideset{}{^\star}\sum}

    \newcommand{\BC}{{\mathbb {C}}} 
     
     \newcommand{\BH}{{\mathbb {H}}}

    \newcommand{\BQ}{{\mathbb {Q}}} \newcommand{\BR}{{\mathbb {R}}}

     \newcommand{\BZ}{{\mathbb {Z}}}

     \newcommand{\GL}{{\mathrm{GL}}}
     \renewcommand{\Im}{{\mathfrak{Im}}}

    \renewcommand{\Re}{{\mathfrak{Re}}}

\def\-{^{-1}}

\def\lp{\left(} 
\def\rp{\right)} 

\def\-{^{-1}}

\newcommand{\delete}[1]{}

     \newcommand{\SL}{{\mathrm{SL}}}
     
     \newcommand{\Sym}{{\mathrm{Sym}}}

    \newcommand{\sstyle}{\scriptstyle}

    \newcommand{\ra}{\rightarrow}

    \theoremstyle{plain}

    \newtheorem{thm}{Theorem}[section] \newtheorem{cor}[thm]{Corollary}
      \newtheorem{prop}[thm]{Proposition}

    \newtheorem {rem}[thm]{Remark}

    \numberwithin{equation}{section}

\begin{document}
\title{Hybrid subconvexity bounds for $L \left(\tfrac{1}{2}, \textnormal{Sym}^2 f \otimes g\right)$} 
\author{R. Holowinsky, R. Munshi, Z. Qi}

\maketitle

\begin{abstract}
Fix an integer $\kappa\geqslant 2$. Let $P$ be prime and let $k> \kappa$ be an even integer. For $f$ a holomorphic cusp form of weight $k$ and full level and $g$ a primitive holomorphic cusp form of weight $2 \kappa$ and level $P$, we prove hybrid subconvexity bounds for $L \left(\tfrac{1}{2}, \textnormal{Sym}^2 f \otimes g\right)$ in the $k$ and $P$ aspects when $P^{\frac {13} {64}  + \delta} < k < P^{\frac 3 8 - \delta}$ for any $0 < \delta < \frac {11} {128}$. These bounds are achieved through a first moment method (with amplification when $P^{\frac {13} {64}} < k \leqslant P^{\frac 4 {13}}$).
\end{abstract}

\section{Introduction}
Subconvexity estimates for Rankin-Selberg $L$-functions have been established in a variety of settings recently with strong motivation coming from equidistribution problems of an arithmetic nature.  In general, for an $L$-function $L(s,\pi)$ associated to an irreducible cuspidal automorphic representation $\pi$ with analytic conductor $Q(s,\pi)$, one hopes to obtain subconvexity estimates of the form $L(s,\pi) \lll Q(s,\pi)^{\frac 1 4-\delta}$ for some $\delta>0$ when $\Re (s)=\tfrac{1}{2}$.  Though the actual value of $\delta$ does not often matter in applications, establishing such a subconvexity bound for some $\delta>0$ is non-trivial and requires careful consideration of the arithmetic/algebraic information associated with $\pi$.  The convexity bound $L(s,\pi) \lll _\varepsilon Q(s,\pi)^{\frac 1 4+\varepsilon}$, on the other hand, follows purely from standard tools in complex analysis.

The resolution of one equidistribution problem related to central values of Rankin-Selberg $L$-functions, the quantum unique ergodicity conjecture of Rudnick and Sarnak \cite{QUE}, has thus far required several techniques from analytic number theory and ergodic theory.  In many cases, however, the conjecture would follow directly from subconvexity estimates for $L\left(\tfrac{1}{2}, \textnormal{Sym}^2 f\right)$ and $L \left(\tfrac{1}{2}, \textnormal{Sym}^2 f \otimes g\right)$.  Here we think of $f$ as a varying modular form and $g$ as a fixed form.

Subconvexity estimates for such $L$-values have proven to be very difficult to establish through current methods and several authors have first given attention to analogous subconvexity problems for Rankin-Selberg $L$-functions in order to possibly better understand the structure behind the symmetric square. For a partial list of related works, see \cite{K-M-V}, \cite{Mich04}, \cite{HM06}, \cite{Ho-Mu}, \cite{HT}, \cite{Blomer12}, \cite{BKY}, \cite{LMY} and the references therein.   For many Rankin-Selberg $L$-functions, it appears as though the arithmetic/analytic structure of the conductor dictates the method of proof that should be adopted to achieve subconvexity.  For example, an amplified moment method is usually required when only one of the forms in the convolution is varying. However, when at least two of the forms are varying, as in the present work, then a moment computation without amplification suffices in certain hybrid ranges. Curiously, the moment method may be avoided all together in cases where the level of the varying form has special structure, for example if the level of the varying form factorizes in a suitable manner \cite{Munshi}.  

In the work of Rizwanur Khan \cite{Khan} on $L \left(\tfrac{1}{2}, \textnormal{Sym}^2 f \otimes g\right)$, a conditional amplifier of long length relative to the conductor was employed in a first moment method in order to establish subconvexity estimates for fixed $f$ and varying $g$ of prime level $P$. Following ideas seen in \cite{Ho-Mu} and \cite{HT} (among others), the work of Khan \cite{Khan} suggests that the number of points of summation for the unamplified first moment of $L \left(\tfrac{1}{2}, \textnormal{Sym}^2 f \otimes g\right)$ is insufficient for application to the subconvexity problem and that one would benefit from increasing the complexity of the $L$-function by allowing $f$ to vary independently with $g$.  

As we demonstrate in this paper, varying the weight $k$ of $f$ along with the level $P$ of $g$ increases the conductor to be of size $Q(\tfrac{1}{2},\textnormal{Sym}^2 f \otimes g)\asymp k^4 P^3$ and allows us to establish hybrid subconvexity bounds for $L \left(\tfrac{1}{2}, \textnormal{Sym}^2 f \otimes g\right)$ in the $k$ and $P$ aspects when ${P^{\frac {13} {64}  + \delta} < k < P^{\frac 3 8 - \delta}}$ for any ${0 < \delta < \frac {11} {128}}$.  Given the above lower bound for $k$ in terms of $P$, this suggests that much more work remains in establishing subconvexity for the case of $P$ fixed and $k$ varying as required in the holomorphic analogue to the quantum unique ergodicity conjecture.  A related situation and hybrid subconvexity bound may be found in (\cite{BKY}, Corollary 1.5) where the authors consider ${L(\tfrac{1}{2}, f \otimes g \otimes h)}$ with all three forms $f,g$ and $h$ varying in weights $k, \ell$ and $k+\ell$ respectively.

\section{Statement of results} \label{sec: Main Theorem}
Fix an integer $\kappa \geqslant 2$ and newform $g_0 \in H^\star_{2\kappa} (P)$ of weight $2\kappa$ and level $P$. Let $f$ be a Hecke eigenform of even weight $k>\kappa $.  Let $L \geqslant 1$ and let $\EuScript{P}$ be a set of primes in the range $[L,2L]$ not dividing the level $P$. Our choice for $\EuScript{P}$ will be such that $|\EuScript{P}|\ggg \frac L { \log L}$.   We will be working with an amplified first moment containing $g_0$
$$
\sum_{g\in H^\star_{2\kappa} (P)} \omega_g^{-1}  |\mathfrak {A}_g|^2 L \left(\tfrac{1}{2}, \textnormal{Sym}^2 f \otimes g\right)
$$
where $\omega_g$ is as in \eqref{omega} and the amplifier is given by
$$
\mathfrak {A}_g :=\sum_{\ell }\alpha_\ell\lambda_g(\ell)
$$
with
\begin{align}
\label{eq: amplifier}
\alpha_\ell :=\begin{cases}\lambda_{g_0}(\ell), &\text{if}\;\ell\in\EuScript{P},\\-1,&\text{if}\;\ell=p^2 \text{ with } p\in\EuScript{P},\\0, &\text{otherwise}.
\end{cases}
\end{align}
When $g=g_0$, the Hecke relation $\lambda_{g_0} (p)^2 - \lambda_{g_0} (p^2) = 1$ yields $ \mathfrak {A}_{g_0} = |\EuScript{P}| $.

Opening the absolute square and using Hecke multiplicativity gives 
\begin{equation} \label{eq: amplified first moment, intro}
\mathop{\sum\sum}_{\ell_1,\ell_2}\alpha_{\ell_1}\overline{\alpha_{\ell_2}}\sum_{\ell_3|(\ell_1,\ell_2)} \sum_{g\in H^\star_{2\kappa} (P)}\omega_g^{-1}\lambda_g(\ell_1\ell_2\ell_3^{-2}) L\left(\tfrac{1}{2}, \textnormal{Sym}^2 f \otimes g \right).
\end{equation}

In \S \ref{sec: proof of Theorem}, we shall prove the following result.

\begin{thm} 
\label{mthm}
Suppose  $k > \kappa \geqslant  2$ are integers, with $ k$ even, $P$ is a prime, and $f$ is a Hecke cusp form of weight $k$ for $\SL(2,\mathbb{Z})$.
Let $\ell\leqslant 16 L^4$ be a positive integer. Then under the assumption
\begin{equation}
\label{eq: assumption on L}
L \leqslant k^{- \frac 2 5} P^{\frac 3 {20}},
\end{equation}
we have for any $\varepsilon > 0$
\begin{align*}
\sum_{g\in H^\star_{2\kappa} (P)} \omega_g^{-1} & \lambda_g(\ell) L \left(\tfrac{1}{2}, \textnormal{Sym}^2 f \otimes g\right) \lll_{\varepsilon, \kappa} \left(\frac{1}{\sqrt{\ell}}+\frac{L^{\frac{22}{7}}k^{\frac{13}{7}}}{P^{\frac{4}{7}}}\right)
(kP)^{\varepsilon}.
\end{align*}
\end{thm}


\bigskip

\begin{rem}
{\rm (1)} The assumptions $\kappa \geqslant 2$ and \eqref{eq: assumption on L} are a result of technical difficulties in the proof. See Remark \ref{rem: Bessel bounds} and \ref{rem: a priori assumption}.

{\rm (2)} Setting $\ell=L=1$, we note that the above bound is the Lindel\"of on average bound when $k \leqslant P^{\frac{4}{13}}$.  Therefore, this is the only range in which amplification is applied.  
\end{rem}

Inserting the above bound into \eqref{eq: amplified first moment, intro} and trivially averaging
over $\ell_1$ and $\ell_2$ we get
\begin{equation*}
\sum_{g\in H^\star_{2\kappa} (P)} \omega_g^{-1}  |\mathfrak {A}_g|^2 L \left(\tfrac{1}{2}, \textnormal{Sym}^2 f \otimes g\right) \lll_{\varepsilon, \kappa} \left( L + \frac {L^{\frac {36} 7} k^{\frac{13} 7}} {P^{\frac 4 7}}\right) (kP)^{\varepsilon}.
\end{equation*}
Using the non-negativity of the central $L$-values and the definition of our amplifier according to \eqref{eq: amplifier} gives
\begin{equation*} 
L\left( \tfrac 1 2, \textnormal{Sym}^2 f\otimes g_0 \right) 
\lll_{\varepsilon, \kappa} \left( \frac P L + L^{\frac {22} 7} k^{\frac {13} 7} P^{\frac 3 7} \right) (kP)^\varepsilon.
\end{equation*}
Finally, setting
\begin{equation}
\label{length-amp}
L = 
\begin{cases}k^{- \frac {13} {29} } P^{\frac 4 {29}}, &\text{ if } P^{\frac {13} {64}} < k \leqslant P^{\frac 4 {13}},\\
1,&\text{ if } P^{\frac 4 {13}} < k < P^{\frac 3 8},
\end{cases}
\end{equation} 
one verifies that the assumption \eqref{eq: assumption on L} on $L$ is satisfied and we therefore obtain the following corollary.
\begin{cor}
For $f$ as above and $g$ a newform of weight $2\kappa$ and level $P$, we have 
\begin{equation} \label{eq: final bound}
L\left( \tfrac 1 2, \textnormal{Sym}^2 f\otimes g \right) 
\lll _{\varepsilon, \kappa}
\begin{cases}
k^{ \frac {13} {29} } P^{\frac {25} {29} } (kP)^{\varepsilon},  & \text { if } P^{\frac {13} {64} } < k \leqslant P^{\frac 4 {13} },\\
 \left( P + k^{\frac {13} 7} P^{\frac 3 7} \right) (kP)^\varepsilon,  & \text{ if } P^{\frac 4 {13} } < k < P^{\frac 3 8}.
\end{cases}
\end{equation}
\end{cor}

\bigskip

\begin{rem}
Note that \eqref{eq: final bound} beats the convexity bound $k P^{\frac 3 4} (kP)^{\varepsilon}$ when $P^{\frac {13} {64}  + \delta} < k < P^{\frac 3 8 - \delta}$ for some $0 < \delta < \frac {11} {128}$.
Putting $\ell = L = 1$ in Theorem \ref{mthm}, we arrive at the following bound by non-negativity $$L\left( \tfrac 1 2, \textnormal{Sym}^2 f\otimes g \right) \lll _{\varepsilon, \kappa} \left( P + k^{\frac {13} 7} P^{\frac 3 7} \right) (kP)^\varepsilon,$$ which is extracted from the second line of \eqref{eq: final bound}. This bound is already able to beat the convexity bound when  $P^{\frac 1 4 + \delta} < k < P^{\frac 3 8 - \delta}$, and therefore amplification is unnecessary (although the bound from amplification, i.e. the first line of \eqref{eq: final bound}, also provides a subconvexity bound on the overlapping range $P^{\frac 1 4 } < k \leqslant P^{\frac 4 {13}}$). Thus, the amplification method extends the range of admissible exponents from below by $\frac 3 {64}$.
\end{rem}

\section{Sketch of hybrid subconvexity in a simplified case}

Let $f$ be a holomorphic cusp form of even weight $k$ and full level and let $g$ be a primitive holomorphic cusp form of even weight $2\kappa$ and prime level $P$.  In order to demonstrate the ideas behind the proofs of our main results, we provide a brief sketch of how one might establish hybrid subconvexity bounds when $\kappa$ is large and fixed.  For notational convenience, we denote the Dirichlet coefficients of $\textnormal{Sym}^2 f$ by $A(n)$ and the coefficients of $g$ by $\lambda(n)$ such that a standard approximate functional equation argument will essentially equate our central $L$-value $L \left(\tfrac{1}{2}, \textnormal{Sym}^2 f \otimes g\right)$ with
$$
D^1_g(Y)+D^2_g(Y):=\sum_{n\lll Y \sqrt{Q}} \frac{A(n)\lambda(n)}{\sqrt{n}}+\varepsilon (\Sym^2 f \otimes g)\sum_{n\lll \sqrt{Q}/Y} \frac{A(n)\lambda(n)}{\sqrt{n}}
$$
for any $Y>0$ with $Q = k^4 P^3 \asymp Q(\tfrac{1}{2},\textnormal{Sym}^2 f \times g)$ and root number $\varepsilon (\Sym^2 f \otimes g)=(-1)^\kappa \sqrt{P}\lambda(P)= \pm 1$.  

Our method will be a normalized first moment average over newforms $g$. Therefore, we wish to achieve a better result than the first moment convexity bound
\begin{equation}\label{FMconvexity}
\sum_{g \in H^\star _{2\kappa}(P)} \omega^{-1}_g L \left(\tfrac{1}{2}, \textnormal{Sym}^2 f \otimes g\right) \lll \frac{Q^{\frac 1 4 + \varepsilon}}{P}
\end{equation}
where $\omega_g$ is as in \eqref{omega}.  As noted in the previous section, one gains from amplification in certain ranges of $k$ relative to $P$, but we omit this component here.

Assume that for our particular choice of $\kappa$ and $P$, the space of newforms $H^\star _{2\kappa}(P)$ spans the space of all forms $S_{2\kappa}(P)$.  Write
$$
S_1(Y)+S_2(Y):=\sum_{g \in H^\star _{2\kappa}(P)} \omega^{-1}_g D^1_g(Y)+\sum_{g \in H^\star _{2\kappa}(P)} \omega^{-1}_g D^2_g(Y)
$$
and consider first $S_1(Y)$.  Applying the Petersson trace formula in the average over $g$ along with standard Bessel function bounds  \eqref{eq: bound for Bessel function} and the Weil bound for Kloosterman sums, one obtains 
\begin{equation*}
\begin{split}
S_1(Y)& = 1+2\pi(-1)^\kappa \sum_{\substack{c\equiv 0 (\mod P)\\ c > 0}} \ \frac{1}{c} \sum_{n\lll  Y\sqrt{Q}} \frac{A(n)}{\sqrt{n}} S(n,1;c) J_{2\kappa-1}\left(\frac{4\pi \sqrt{n}}{c}\right)\\
& \lll _\kappa  1+ \sum_{\substack{c\equiv 0 (\mod P)\\ c > 0}} \ \frac{1}{\sqrt{c}} \sum_{n\lll  Y\sqrt{Q}} \frac{|A(n)|}{\sqrt{n}}\left(\frac{\sqrt{n}}{c}\right)^{2\kappa-1}.
\end{split}
\end{equation*}
The Deligne bound for the coefficients of holomorphic forms then gives
\begin{equation}\label{S1sketch}
S_1(Y) \lll _{\varepsilon, \kappa} \max \left\{1,\sqrt{P} \left(\frac{Y \sqrt{Q}}{P^2}\right)^\kappa \right\}Q^{\varepsilon}.
\end{equation}
Note that the condition $1<\frac{Q^{\frac 1 4}}{P}$ or equivalently $P^{\frac 1 4}<k$ is necessary for the bound \eqref{S1sketch} above to be better than the first moment convexity bound \eqref{FMconvexity}.
\begin{rem}
We shall see that such basic analysis of $S_1(Y)$, when $\kappa$ is large, is sufficient for establishing at least some hybrid subconvexity bound due to the congruence condition in the sum over $c$.  In \S \ref{sec: proof of Theorem}, we improve on the above bound \eqref{S1sketch} in order to establish our main results in \S \ref{sec: Main Theorem}.
\end{rem}

Now for $S_2(Y)$, an application of Petersson's trace formula shows that $S_2(Y)$ is essentially equal to
\begin{equation*}
\begin{split}
& \sum_{\substack{c\equiv 0 (\mod P)\\ c > 0}} \frac{\sqrt{P}}{c} \sum_{n\lll  \sqrt{Q}/Y}\frac{A(n)}{\sqrt{n}} S(nP,1;c) J_{2\kappa-1}\left(\frac{4\pi \sqrt{nP}}{c}\right)\\
= &  \sum_{\substack{(c,P)=1\\ c > 0}} \frac{1}{c\sqrt{P}} \sum_{n\lll  \sqrt{Q}/Y}\frac{A(n)}{\sqrt{n}} S(nP,1;cP) J_{2\kappa-1}\left(\frac{4\pi \sqrt{n}}{c\sqrt{P}}\right)\\
= &  - \sum_{\substack{(c,P)=1\\ c > 0}} \frac{1}{c\sqrt{P}} \sum_{n\lll  \sqrt{Q}/Y}\frac{A(n)}{\sqrt{n}} S(n,\overline{P};c) J_{2\kappa-1}\left(\frac{4\pi \sqrt{n}}{c\sqrt{P}}\right).
\end{split}
\end{equation*}
Here we pulled out the $P$ divisor in the original $c$-sum and used basic properties of the Kloosterman sums, i.e. 
$$
S(nP,1;cP)=\begin{cases}
0, & \textnormal{ if } P|c,\\
S(0,\overline{c};P)S(n,\overline{P};c)=-S(n,\overline{P};c), & \textnormal{ if } (c,P)=1.
\end{cases}
$$
Focusing on the transition range of the Bessel function ($\sqrt{n} \asymp c \sqrt{P}$) for the remainder of the sketch and opening the Kloosterman sum, we see that we must analyze a smoothed version of
$$
\sum_{\substack{(c,P)=1\\ c \asymp Q^{\frac 1 4}/\sqrt{YP}}} \frac{1}{c\sqrt{P}}\ \sideset{}{^\star}\sum_{a(\mod c)} e\left(\frac{a \overline{P}}{c}\right)\sum_{n\lll \sqrt{Q}/Y}\frac{A(n)}{\sqrt{n}} e\left(\frac{n\overline{a}}{c}\right).
$$
An application of Vorono\"{i} summation in $n$ leads to sums of the form
$$
\sum_{\substack{(c,P)=1\\ c \asymp Q^{\frac 1 4}/\sqrt{YP}}} \frac{1}{c\sqrt{P}}\ \sideset{}{^\star}\sum_{a(\mod c)} e\left(\frac{a \overline{P}}{c}\right)\sum_{n\lll \frac{c^3k^2}{\sqrt{Q}/Y}}\frac{A(n)}{\sqrt{n}} \frac{S(a,n;c)}{\sqrt{c}}
$$
i.e. we obtain a new $n$-sum of length ``conductor divided by the original length of summation'' with a summand ``dual'' to the previous summand.  Summing over $a$, one sees that we must consider
$$
\sum_{\substack{(c,P)=1\\ c \asymp Q^{\frac 1 4}/\sqrt{YP}}} \frac{1}{\sqrt{cP}}\sum_{n\lll \frac{ k^2 Q^{\frac 1 4}}{\sqrt{YP^3}}}\frac{A(n)}{\sqrt{n}} e\left(-\frac{nP}{c}\right).
$$
Trivially bounding the sum over $n$ using Deligne's bound, the above is bounded by
\begin{equation}\label{S2sketch}
\sum_{\substack{(c,P)=1\\ c \asymp Q^{\frac 1 4}/\sqrt{YP}}} \frac{1}{\sqrt{cP}}\left(\frac{k Q^{\frac 1 8}}{Y^{\frac 1 4} P^{\frac 3 4}}\right)^{1+\varepsilon} \lll \frac{Q^{\frac 1 4+\varepsilon}}{P}\left(\frac{k}{\sqrt{YP}}\right).
\end{equation}
This bound is better than the first moment convexity bound \eqref{FMconvexity} when $Y>\frac{k^2}{P}$.

We now combine the bounds in \eqref{S1sketch} and \eqref{S2sketch}.  Assume first that 
$$
\sqrt{P} \left(\frac{Y \sqrt{Q}}{P^2}\right)^\kappa \leqslant 1<\frac{Q^{\frac 1 4}}{P}.
$$  
Equating the two bounds $1=\frac{Q^{\frac 1 4}}{P}\left(\frac{k}{\sqrt{YP}}\right)$ we get
$Y=\frac{ k^2 \sqrt{Q}}{P^3}= k^4 P^{-\frac 3 2} $.  Such a choice of $Y$ satisfies our assumption when $P^{\frac 1 4}< k \leqslant P^{\frac{1}{3}-\frac{1}{12\kappa}}$.  Now assume that 
$$
1\leqslant \sqrt{P} \left(\frac{Y \sqrt{Q}}{P^2}\right)^\kappa<\frac{Q^{\frac 1 4}}{P}.
$$
Equating the two bounds $\sqrt{P} \left(\frac{Y \sqrt{Q}}{P^2}\right)^\kappa=\frac{Q^{\frac 1 4}}{P}\left(\frac{k}{\sqrt{YP}}\right)$ we get 
$
Y= k^{\frac{4-4\kappa}{2\kappa+1}} P^{\frac{2\kappa-5}{4\kappa+2}}.
$
Such a choice of $Y$ satisfies our assumption when $P^{\frac{1}{3}-\frac{1}{12\kappa}} \leqslant k<P^{\frac{3(2\kappa-1)}{4(4\kappa-1)}}$.  

Therefore, one establishes hybrid subconvexity bounds for all $\kappa\geqslant 2$ with the range of $k$ relative to $P$ tending to $P^{\frac 1 4} < k < P^{\frac 3 8}$ as $\kappa \longrightarrow \infty$.


\section{Preliminaries}

\subsection{Holomorphic cusp forms} For a positive integer $N$ and an even positive integer $k$, the space $S_k(N)$ of cusp forms of weight $k$ for the Hecke congruence group $\Gamma_0 (N)$ is a finite-dimensional Hilbert space with respect to the Petersson inner product
$$\langle f_1, f_2\rangle := \int_{\Gamma_0(N) \backslash \BH } f_1(z) \overline {f_2(z)} y^{k-2} dx d y, \hskip 10 pt f_1, f_2 \in S_k(N),$$ where $\BH$ denotes the upper half-plane. 
Every $f \in S_k(N)$ has a Fourier series expansion
$$f(z) = \sum_{n = 1}^\infty \psi_f (n) n^{ \frac {k-1} 2} e(nz).$$

For $n \geqslant  1$, define the Hecke operator $T_N(n)$ by
$$(T_N(n) f) (z) := \frac 1 {\sqrt n} \sum_{\sstyle ad = n\atop \sstyle (a, N) = 1} \left( \frac a d \right)^{\frac k 2} \sum_{b (\mod d)} f\left( \frac {az + b} d \right).$$

Let $H^\star_k(N)$ be the orthogonal set of Hecke-normalized (i.e., $\psi_f(1) = 1$) newforms $f$ in $S_k(N)$. Every $f \in H^\star_k(N)$ is an eigenfunction of all Hecke operators $T_N(n)$; let $\lambda_f(n)$ be its eigenvalue of $T_N(n)$. We have $\psi_f(n) = \lambda_f(n)$ for all $n \geqslant  1$. The Hecke eigenvalues are multiplicative, i.e., for any $m , n \geqslant  1$
\begin{equation}
\label{eq: Hecke eigenvalues multiplicative}
\lambda_f(m) \lambda_f(n) = \sum_{\sstyle d  |  (m ,n) \atop \sstyle (d, N) = 1} \lambda_f (mn/d^2).
\end{equation} 
In particular, (\ref{eq: Hecke eigenvalues multiplicative}) becomes completely multiplicative when $n  |  N$ 
\begin{equation}
\label{eq: Hecke eigenvalues multiplicative 2}
\lambda_f(m) \lambda_f(n) = \lambda_f (mn).
\end{equation}
For any $f \in H^\star_k(N)$, we have Deligne's bound
\begin{equation*} 
|\lambda_f(n)| \leq  \tau (n),
\end{equation*}
and when $N$ is squarefree, it is known that (\cite[(2.24)]{ILS})
\begin{equation*} 
\lambda_f(n)^2 =  \frac 1 n, \hskip 10pt \text { if } n |  N.
\end{equation*}

\subsection{Automorphic $L$-functions}

In this section some preliminary results on automorphic $L$-functions are given. We shall particularly focus on the Rankin-Selberg $L$-function $L(s, \Sym^2 f\otimes g )$ for  $f \in S _k(1)$ and $g \in  H^\star _{2\kappa}(N)$
with $N$ squarefree, $k$ an even positive integer and $\kappa $ a positive integer. A brief calculation of the $\gamma$-factor and the $\varepsilon$-factor of $\Sym^2 f\otimes g$ will be given in \S \ref{sec: Computations of gamma and epsilon}.

\subsubsection{ } For $f \in H^\star_k(N)$ the Hecke $L$-function is defined by
\begin{equation*}
L(s, f) = \sum_{n = 1}^\infty \lambda_f(n) n^{-s}.
\end{equation*}
This has an Euler product $L(s, f) = \prod_{p} L_p(s, f)$ with local factors
$$L_p(s, f) = (1 - \lambda_f(p) p^{-s} + \chiup_0 (p)^{-2s})\-, $$
where $\chiup_0 $ is the principal character modulus $N$.
The gamma factor is
\begin{equation*}
\gamma (s, f) = \pi^{-s} \Gamma \left( \frac { s} 2 + \frac {k - 1} 4 \right)  \Gamma \left( \frac { s } 2  + \frac {k + 1} 4 \right).
\end{equation*}
The complete product $\Lambda(s, f) = N^{\frac s 2}\gamma (s, f) L(s, f)$ is entire and satisfies the functional equation
\begin{equation*}
\Lambda (s, f) = \varepsilon_f \Lambda (1-s, f),
\end{equation*}
with root number $\varepsilon_f = i^k \eta_f = \pm 1,$ where $\eta_f$ is the eigenvalue of the Atkin-Lehner involution $W_N$. If $N$ is squarefree, then $\varepsilon(f) = i^k \mu(N) \lambda_f (N) \sqrt N.$

For $p \notdivides N$, the local factors $L_p(s, f)$ factor further as
\begin{equation*}
L_p(s, f) = \left( 1 - \alpha_f(p) p^{-s} \right)\- \left( 1 - \beta_f(p) p^{-s} \right)\-,
\end{equation*}
where $\alpha_f(p)$, $\beta_f(p)$ are complex numbers with $\alpha_f(p) + \beta_f(p) = \lambda_f (p)$ and $\alpha_f(p) \beta_f(p) = 1$.

\subsubsection{ } For $f \in S_k(1)$ the symmetric square $L$-function is defined by
\begin{equation*}
L(s, \Sym^2 f) = \zeta (2 s) Z(s, f),
\end{equation*}
where $\zeta (s)$ denotes the Riemann zeta function and $Z(s, f)$ is defined by $$Z(s, f) = \sum_{n = 1}^\infty \lambda_f(n^2) n^{-s}.$$
This has Euler product $L(s, \Sym^2 f) = \prod_p L_p(s, \Sym^2 f)$ with
\begin{equation*}
L_p(s, \Sym^2 f) = (1 - \alpha_f^2 (p) p^{-s})\- (1 - p^{-s})\- (1 - \beta_f^2 (p) p^{-s})\-.
\end{equation*}
The gamma factor is
\begin{equation*}
\gamma (s, \Sym^2 f) = \pi^{- \frac {3s} 2 } \Gamma \left( \frac {s+ 1} 2 \right)  \Gamma \left( \frac {s + k - 1} 2 \right)  \Gamma \left( \frac {s+ k} 2 \right) 
\end{equation*}
The complete product $\Lambda(s, \Sym^2 f) = \gamma (s, \Sym^2 f) L(s, \Sym^2 f)$ is entire and it satisfies the functional equation
\begin{equation*}
\Lambda (s, \Sym^2 f) = \Lambda (1 - s, \Sym^2 f).
\end{equation*}

We have the convexity bound
\begin{equation} \label{eq: L(s, Sym2 f) convexity bound}
L(\sigma + i t, \Sym^2 f) \lll_\varepsilon \left( (|t| + 1) (|t| + k)^2 \right)^{\frac { 1 - \sigma} 2 + \varepsilon}, \hskip 10pt 0 \leqslant \sigma \leqslant 1,
\end{equation}
where the implied constant depends only on $\varepsilon > 0$.
Moreover, it is known that \cite{HL}
\begin{equation} \label{eq: L(1, Sym2 f)}
k^{- \varepsilon} \lll_\varepsilon L(1, \Sym^2 f) \lll_\varepsilon k^{ \varepsilon}.
\end{equation}

According to \cite{GJ}, $L(s, \Sym^2 f)$ is also an $L$-function $L(s, F)$ of some automorphic representation $F$ of $\GL(3, \BZ)$, and the normalized Fourier coefficients are given by
\begin{equation*}
A_F(1, n) = A_F(n, 1) = \sum_{m \ell^2 = n} \lambda_f(m^2),
\end{equation*}
\begin{equation*}
A_F(m, n) = A_F(- m, n) =  A_F( m, - n) =  A_F(- m, - n),
\end{equation*}
\begin{equation*}
 A_F(m, 0) = A_F(0, n) = 0,
\end{equation*}
and the Hecke relations (\cite[(7.7)]{MS})
\begin{equation} \label{eq: Hecke relation}
A_F(m, n) = \sum_{d |  (m, n)} \mu(d) A_F(m/d, 1) A_F(1, n/d).
\end{equation}
We have
\begin{equation*}
L(s, \Sym^2 f) = L(s, F) = \sum_{n=1}^\infty A_F(1, n) n^{-s}.
\end{equation*}

\subsubsection{ }
Let  $f \in S _k(1)$ and $g \in  H^\star _{2\kappa}(N)$
with $N$ squarefree, $k$ an even positive integer and $\kappa $ a positive integer. We define the Rankin-Selberg $L$-function
\begin{equation*} 
L(s, \Sym^2 f \otimes g) = \sum_{m = 1}^\infty  \sum_{n = 1}^\infty A_F(m, n) \lambda_g(n) (m^2 n)^{-s}.
\end{equation*}
This has Euler product $L(s, \Sym^2 f \otimes g) = \prod_p L_p(s, \Sym^2 f \otimes g)$ with local factor
\begin{equation*}
\begin{split}
& \quad \ L_p(s, \Sym^2 f \otimes g) \\ 
& = (1 - \alpha^2_f(p)\alpha_g(p) p^{-s})\- (1 - \alpha_g(p) p^{-s})\- (1 - \beta^2_f(p)\alpha_g(p) p^{-s})\- \\
& \hskip 15pt (1 - \alpha^2_f(p)\beta _g(p) p^{-s})\- (1 - \beta _g(p) p^{-s})\- (1 - \beta^2_f(p)\beta _g(p) p^{-s})\-
\end{split}
\end{equation*}
if $p \notdivides  N$, and
\begin{equation*}
\begin{split}
& \ L_p(s, \Sym^2 f \otimes g) \\
= &(1 - \alpha^2_f(p)\lambda_g(p) p^{-s})\- (1 - \lambda_g(p) p^{-s})\- (1 - \beta^2_f(p)\lambda_g(p) p^{-s})\-
\end{split}
\end{equation*}
if $p  |  N$.
The gamma factor is
\begin{equation} \label{eq: gamma factor for Sym2 f otimes g}
\begin{split}
& \ \gamma(s, \Sym^2 f \otimes g) \\
=& (2\pi)^{-3 s} \Gamma \left(s + \kappa - \frac { 1} 2\right) \Gamma \left( s + k + \kappa - \frac { 3} 2 \right) 
\Gamma \left( s + \left| k - \kappa - \frac 1 2 \right| \right).
\end{split}
\end{equation}
The complete product $\Lambda(s, \Sym^2 f \otimes g) = N^{\frac {3s} 2 } \gamma (s, \Sym^2 f \otimes g) L(s, \Sym^2 f \otimes g)$ is entire and it satisfies the functional equation
\begin{equation*}
\Lambda (s, \Sym^2 f \otimes g) = \varepsilon (\Sym^2 f \otimes g) \Lambda (1 - s, \Sym^2 f \otimes g),
\end{equation*}
with root number $\varepsilon (\Sym^2 f \otimes g) = \pm 1$ given by
\begin{equation}  \label{eq: epsilon factor of Sym2 f otimes g}
\varepsilon (\Sym^2 f \otimes g) = \left \{ \begin{split}
& (- 1)^{\kappa + 1} \mu (N) \lambda_g(N) \sqrt N, \hskip 10pt & \text { if } k > \kappa,\\
& (- 1)^{\kappa } \mu (N) \lambda_g(N) \sqrt N, \hskip 10pt & \text { if } k \leqslant \kappa.
\end{split}\right.
\end{equation}

It follows from the Hecke relations (\ref{eq: Hecke relation}) and (\ref{eq: Hecke eigenvalues multiplicative 2}) that
\begin{equation} \label{eq: L(s, Sym2 f otimes g 2)}
\begin{split}
& \quad \ L(s, \Sym^2 f \otimes g) \\
& = \sum_{d = 1}^\infty \mu(d) d^{-3 s} \sum_{m = 1}^\infty A_F(m, 1) m^{-2 s} \sum_{n = 1}^\infty A_F(1, n) \lambda_g(d n) n^{-s}\\
& = L(2s, \Sym^2 f ) \sum_{d = 1}^\infty \mu(d) d^{-3 s} \sum_{n = 1}^\infty A_F(1, n) \lambda_g(d n) n^{-s}\\
& =  L(2s, \Sym^2 f ) L_N (3s, g)\- \sum_{(d, N) = 1} \mu(d) d^{-3 s} \sum_{n = 1}^\infty A_F(1, n) \lambda_g(d n) n^{-s},
\end{split}
\end{equation}
where $L_N (s, g)$ denotes the finite Euler product
\begin{equation*}
L_N ( s, g) = \prod_{p | N} \left(1 - \lambda_g(p) p^{-s} \right)\-.
\end{equation*}

\subsubsection{Computation of the $\gamma$-factor and the $\varepsilon$-factor of $\Sym^2 f \otimes g$}\label{sec: Computations of gamma and epsilon}

Let  $\psi_{\infty}$ be the standard additive character on $\BR$, namely $\psi_\infty (x) = e(x)$, and $\psi_p$ be a normalized unramified additive character of $\BQ_p$ for each prime $p$.

At the real place, the local component $\pi_{f, \infty}$, respectively $\pi_{g, \infty}$, is the discrete series representation of $\GL(2, \BR)$ with weight $k$, respectively $2 \kappa$.

The Weil group The Weil group $W_{\BR} $ of $\BR$ is realized as $ \BC^\times \cup j \BC^\times$ satisfying $j^2 = - 1 \in \BC^\times$ and $j  z j\- = \overline z$ for $z \in \BC^\times$. Under the local Langlands correspondence the discrete series with weight $k$,  $k \geqslant  2 $, corresponds to the two dimensional representation $\rho_k$ of $W_{\BR}$ given by
\begin{equation*}
\rho_k (r e^{i \theta}) = \begin{pmatrix}
e^{i (k - 1) \theta} &  \\
 & e^{- i (k-1) \theta}
\end{pmatrix}
\hskip 10pt \rho_k (j) = \begin{pmatrix}
 & (-1)^{k-1} \\
1 &
\end{pmatrix}
\end{equation*}
The $\gamma$-factor $ \gamma (s, \rho_k) = \Gamma_{\BC} \left( s + \frac {k - 1} 2 \right)$ and the $\varepsilon$-factor $ \varepsilon (\rho_k, \psi_\infty) = i^k$. See \cite{Knapp}.
With some matrix calculations we have
\begin{equation*}
\Sym^2 (\rho_k) \otimes \rho_{2\kappa} = \rho_{2\kappa} \oplus \rho_{|2 k - 2\kappa - 1| + 1} \oplus \rho_{2k + 2\kappa - 2}.
\end{equation*}
This implies the formula (\ref{eq: gamma factor for Sym2 f otimes g}) for the $\gamma$-factor $\gamma (s, \Sym^2 f \otimes g)$, and the $\varepsilon$-factor at $\infty$ is
\begin{equation}\label{eq: epsilon factor at infty}
\varepsilon_{\infty} (\Sym^2 f \otimes g, \psi_\infty) = \left \{ \begin{split}
&(-1)^{\kappa + 1}, \hskip 10pt & \text { if } k > \kappa,\\
& (-1)^{\kappa }, \hskip 10pt & \text { if } k \leqslant \kappa.
\end{split} \right.
\end{equation}

For any prime $p $, the local component $\pi_{\Sym^2 f, p}$ is an unramified principle series representation of $\GL(3, \BQ_p)$ with trivial central character, so
\begin{equation}\label{eq: epsilon factor at p}
\varepsilon_p  (\Sym^2 f \otimes g, \psi_p) = \varepsilon_p  (g, \psi_p)^3 = \left \{ \begin{split}
& - \lambda_g(p) \sqrt p, \hskip 10pt \text { if } p | N,\\
& 1, \hskip 60pt \text { if } p \notdivides N.
\end{split} \right.
\end{equation}
Multiplying (\ref{eq: epsilon factor at infty}) and (\ref{eq: epsilon factor at p}) yields the formula (\ref{eq: epsilon factor of Sym2 f otimes g}) for the $\varepsilon$-factor $\varepsilon (\Sym^2 f \otimes g)$.

\delete{
(Indeed, $\Sym^2 (\rho_k) = \rho_0^- \oplus \rho_{2k - 1}$ with $\rho_0^-$ defined by $\rho_0^- (z) = 1$ and $\rho_0^- (j) = - 1$. For instance, \begin{equation*}
\begin{split}
\Sym^2 (\rho_k) (r e^{i \theta})&  = \begin{pmatrix}
e^{i (2k - 2) \theta} & & \\
 & 1 & \\
 & & e^{- i (2k-2) \theta}
\end{pmatrix}, \\
\Sym^2 (\rho_k) (j) & = \begin{pmatrix}
 & & (-1)^{2k-2} \\
 & (-1)^{k - 1} & \\
 1 & & 
\end{pmatrix},
\end{split}
\end{equation*}
and
\begin{equation*}
\begin{split}
(\rho_{2k-1} \otimes \rho_{2\kappa}) (r e^{i \theta})&  = \begin{pmatrix}
e^{i (2k + 2\kappa - 3) \theta} & & & \\
 & e^{i (2k - 2\kappa - 1) \theta} & & \\
 & & e^{- i (2k- 2\kappa - 1) \theta} & \\
 & & & e^{- i (2k + 2\kappa - 3) \theta}
\end{pmatrix}, \\
(\rho_{2k-1} \otimes \rho_{2\kappa}) (j) & = \begin{pmatrix}
& & & (-1)^{2k + 2\kappa - 3} \\
& & (-1)^{2k - 2} & \\
& (-1)^{2\kappa - 1} & & \\
1 & & &
\end{pmatrix}.
\end{split}
\end{equation*}
)
}

\subsection{Approximate functional equation}

In view of (\ref{eq: L(s, Sym2 f otimes g 2)}), we have the following approximate functional equation (see 
\cite[Theorem 5.3, Proposition 5.4]{IK})
\begin{equation}
\label{eq: Approximate functional equation}
\begin{split}
& \quad \  L \left(\tfrac 1 2, \Sym^2 f \otimes g \right) \\
& =  \sum_{(d, N) = 1} \sum_{n = 1}^\infty \frac { \mu(d) A_F(1, n)  \lambda_g(d n)} { \sqrt {d^{ 3 } n }} V\left( \frac {d^3 n} {Y N^{\frac 3 2}}\right) \\
& \quad \ + \varepsilon (\Sym^2 f \otimes g) \sum_{(d, N) = 1} \sum_{n = 1}^\infty \frac { \mu(d) A_F(1, n)  \lambda_g(d n)} { \sqrt {d^{ 3 } n }} V\left( \frac {d^3 n Y} {N^{\frac 3 2}}\right),
\end{split}
\end{equation}
where $V(y)$ is a smooth function on $\BR_+$ defined by
\begin{equation*}
\begin{split}
\frac 1  {2 \pi i} \int_{(3)} \left( \cos \lp \frac {\pi u} {4A} \rp \right)^{- 24 A} \frac {\gamma \left(\frac 1 2 + u, \Sym^2 f \otimes g \right)}  {\gamma \left(\frac 1 2, \Sym^2 f \otimes g \right)} \frac { L(1 + 2u, \Sym^2 f) } { L_N \left(\frac 3 2 + 3 u, g \right)} y^{ - u} \frac { du } u,
\end{split}
\end{equation*}
with $A$ a positive integer.
For $y$ large we shift the contour of integration in $V(y)$ to $\Re\ u = A$, and for $y$ small we left shift the contour  to $\Re\ u = - \sigma$ with $0 < \sigma < \frac 1 2$, passing through the pole $u = 0$ with residue $L(1, \Sym^2 f) L_N \left(\frac 3 2, g \right)\-$. Then by Stirling's formula and the convexity bound (\ref{eq: L(s, Sym2 f) convexity bound}) for $L(s, \Sym^2 f)$ we derive
\begin{equation}\label{eq: V(y)}
y^j V^{(j)} (y) \lll_{\varepsilon, j, A, \kappa} (k N)^\varepsilon \left( 1 + \frac y {k^2} \right)^{- A},
\end{equation}
and the asymptotic equation
\begin{equation}\label{eq: V(y) 2}
y^j V^{(j)} (y) = L(1, \Sym^2 f) L_N \left(\tfrac 3 2, g \right)\- \delta (0, j) + O_{\varepsilon, \sigma, j, \kappa} \left( (kN)^\varepsilon y ^{\sigma} \right).
\end{equation}
Choosing $\sigma$ sufficiently small, we have on the range $ y < k^{2 + \varepsilon}$
\begin{equation}\label{eq: V(y) y < k2}
y^j V^{(j)} (y) \lll_{\varepsilon, j, \kappa} (k N)^\varepsilon.
\end{equation}


\subsection{Petersson trace formula}

Let $\mathscr B _k (N)$ be an orthogonal basis of $S_k(N)$. For any $n, m \geqslant  1$ define
\begin{equation*} 
\Delta_{k, N} (m, n) := \sum_{f \in \mathscr B_k (N)} \omega_f\- \psi_f(m) \psi_f(n).
\end{equation*}
This is basis independent. Here the weight $\omega_f$ is defined by 
\begin{equation} \label{omega}
\omega_f := \frac {(4 \pi)^{k - 1}} {\Gamma (k - 1)} \langle f, f \rangle.
\end{equation}
If $f \in H^\star_k (N)$, then 
$$\omega_f = \frac {(k-1) N} {2 \pi^2} L(1, \Sym^2 f),$$ 
and by (\ref{eq: L(1, Sym2 f)}) we have
\begin{equation} \label{eq: bound of omega f}
(kN)^{1 - \varepsilon } \lll_\varepsilon  \omega_f \lll_\varepsilon (kN)^{1 + \varepsilon}.
\end{equation}

Moreover, it follows from \cite[(2.48, 2.72)]{ILS}, the bound (\ref{eq: L(1, Sym2 f)}) of $L(1, \Sym^2 f)$ and Deligne's bound that 
\begin{equation} \label{eq: trivial bound for DeltakN}
\Delta_{k, N} (m, n) \lll_\varepsilon ( N k n m)^\varepsilon.
\end{equation}

We have the following formula of Petersson. 
\begin{equation} \label{eq: Petersson's formula}
\Delta_{k, N} (m, n) = \delta(m, n) + 2 \pi i^{-k} \sum_{\sstyle c > 0 \atop \sstyle c \equiv 0 (\mod N)}  \frac { S(m, n; c)} c J_{k-1} \left( \frac {4 \pi \sqrt {m n}} c \right),
\end{equation}
where $J_{k-1}$ is the Bessel function of the first kind with order $k-1$.

Define 
\begin{equation*} 
\Delta^{\star}_{k, N} (m, n) := \sum_{f \in H^\star_k (N)} \omega_f\- \lambda_f(m) \lambda_f(n).
\end{equation*}
According to \cite[Proposition 2.8]{ILS}, under the assumptions that $N$ be squarefree, $(m, N) = 1$ and $(n, N) | N^2$ we have
\begin{equation} \label{eq: Petersson for Delta star, original}
\Delta^\star _{k, N} (m, n) =  \sum_{LR = N} \frac {\mu(L) } {L \nu((n, L))} \sum_{\ell  |   L^\infty} \ell\- \Delta_{k, R} \left(m \ell^2, n \right),
\end{equation}
with
\begin{equation*} 
\nu (n) =  n \prod_{p  |   n} \left( 1 + \frac 1 p \right). 
\end{equation*}
For our purpose, a more general variant of (\ref{eq: Petersson for Delta star, original}) given in \cite[Lemma 2.4]{HT}  is more convenient. Suppose $N$ is squarefree and $(m, N) = 1 $, then
\begin{equation}\label{eq: Petersson for Delta star}
\Delta^\star _{k, N} (m, n) = \sum_{LR = N} \frac {\mu(L) } {L \nu((n, L))} \sum_{\ell  |   L^\infty} \ell\- \sum_{\ell_1^2  |   (n, \ell_1 L)} \mu(\ell_1) \ell_1 \Delta_{k, R} \left(m \ell^2, n {\ell_1^{- 2}} \right).
\end{equation}

\subsection{Bessel functions}
Let $\nu$ be a positive integer. If $x \lll 1$, the Taylor series expansion yields
\begin{equation}
\label{eq: bound for Bessel function}
x^j J_{\nu }^{(j)} (x) \lll_{\nu, j} x^{\nu}, \hskip 10 pt j \geqslant  0.
\end{equation}
We have (see \cite[\S 16.12, \S 17.5]{Whittaker-Watson})
\begin{equation} \label{eq: Bessel function and Whittaker function}
J_{\nu}(x) = \frac {1} {\sqrt {2 \pi x}} \left( e^{i x} W_{\nu, +}(x) + e^{- i x} W_{ \nu, -} (x) \right),
\end{equation}
where
\begin{equation*}
\begin{split}
W_{\nu, \pm}(x) & = e^{\mp i  (x + \frac 1 2\lp \nu + \frac 1 2 \rp \pi i )} W_{0, \nu} (\mp 2 i x) \\
& = \frac {e^{\mp \frac 1 2\lp \nu + \frac 1 2 \rp \pi i  }} {\Gamma\left(\nu + \frac 1 2 \right) }  \int_0^\infty e^{-y} \left( y \left( 1 \pm i \frac {y } {2x} \right) \right) ^{\nu - \frac 1 2 } d y.
\end{split}
\end{equation*}

For $x \ggg 1$, the asymptotic expansion for the Whittaker functions (\cite[\S 16.3]{Whittaker-Watson} or \cite[9.227]{GR}) and their recursion formula (\cite[9.234.3]{GR}) provide the bound
\begin{equation}
\label{eq: bound of Whittaker functions}
x^j W_{\nu, \pm}^{(j)} (x) \lll_{\nu, j} 1.
\end{equation}

\begin{rem}\label{rem: Bessel bounds}
The authors would like to point out a common mistake presented in literatures. \eqref{eq: bound of Whittaker functions} does {\rm not} hold for small $x \lll 1$. Actually, we have
$$J_\nu (x) + i Y_\nu (x) = H_{\nu}^{(1)} (x) = \sqrt {\frac {2} {\pi x}} e^{i x} W_{\nu, +}(x),$$
where $Y_\nu$ is the Bessel function of the second kind and $H^{(1)}_{\nu}$ is the Hankel function. The behaviour of $Y_\nu$ near zero is given by \cite[3.52 (3)]{Watson} as follows
\begin{equation*}
\begin{split}
Y_{\nu} (x) = & - \frac 1 {\pi} \sum_{n=0}^{\nu - 1} \frac {(\nu - n - 1)!} {n!} \lp \tfrac 12 x \rp^{- \nu - 2 n} \\
& + \frac 1 \pi \sum_{n=0}^\infty \frac {(-1)^n \lp \tfrac 12 x \rp^{\nu + 2 n}} {n! (\nu + n)!} \lp 2 \log \lp \tfrac 12 x \rp - \psi (n+1) - \psi (\nu + n + 1) \rp,
\end{split}
\end{equation*} 
with $\psi (1) = - \gamma$, $\psi (n+1) = 1 + \frac 1 2 + ... + \frac 1 m - \gamma$, where $\gamma $ denotes Euler's constant. In particular, $W_{\nu, +} (x)$ tends to infinity with growth rate $x^{\frac 1 2 - \nu}$ as $x$ tends to zero.

Thus, for small $x$, one has to apply the bound \eqref{eq: bound for Bessel function} of $J_\nu (x)$, as done for instance in \S \ref{sec: theta sqrt X < 1}. However, $\nu \geqslant 3$, or $\kappa \geqslant 2$ in our context, is required to guarantee convergence of certain series as shown in {\rm (\ref{eq: c-sum requiring kappa > 1, 1}, \ref{eq: c-sum requiring kappa > 1, 2})}.
\end{rem}

\subsection{Vorono\"i formula} \label{sec: Voronoi}
For a smooth compactly-supported function $\psi(y)$ on $\BR_+$, define the Mellin transform by
\begin{equation*} 
\widetilde {\psi} (s) := \int_0^\infty \psi(y) y^{s } \frac { d y } y.
\end{equation*}
Let $\eta \in \{0, 1 \}$. For $\sigma > - 1$ define
\begin{equation*} 
\Psi_\eta (x) := \frac 1 { 2 \pi i} \int_{(\sigma)} (\pi^3 x)^{-s} G_{\eta} (s) \widetilde  {\psi} (- s) d s,
\end{equation*}
where
\begin{equation}
\label{eq: definition G eta (s)}
G_{\eta} (s) := \frac 
{\Gamma \left( \frac {1 + s + 1 - \eta} 2 \right) \Gamma \left( \frac { 1 + s + (k-1) + \eta} 2 \right) \Gamma \left( \frac { 1 + s + k - \eta} 2 \right)} 
{\Gamma \left( \frac { - s + 1 - \eta} 2 \right) \Gamma \left( \frac { - s + (k-1) + \eta} 2 \right)  \Gamma \left( \frac { - s + k - \eta} 2 \right)}.
\end{equation}
Then define
\begin{equation*}
\begin{split}
\Psi_+ (x) & := \frac 1 { 2 \pi^{\frac 3 2}} \left( \Psi_0 (x) - i \Psi_1 (x) \right),\\
\Psi_- (x) & := \frac 1 { 2 \pi^{\frac 3 2}} \left( \Psi_0 (x) + i \Psi_1 (x) \right).
\end{split}
\end{equation*}
We have the following Vorono\"i formula (\cite[Theorem 1.18]{MS}).

\begin{prop}
Let $\psi (y)$ be a smooth function compactly supported on $\BR_+$. Let $a, \overline a, c \in \BZ $ with $c \neq 0$, $(a, c) = 1$ and $a \overline a \equiv 1 (\mod c)$. Then we have
\begin{equation*}
\begin{split}
\sum_{n = 1}^\infty A_F(1, n) e \left( \frac {n \overline a} c \right) \psi (n) = & c \sum_{n_1 | c} \sum_{n_2 = 1}^\infty \frac {A_F(n_2, n_1)} {n_1 n_2} S\left(a, n_2;  c / {n_1} \right) \Psi_+\left( \frac {n_2 n_1^2} {c^3} \right)\\
& +  c \sum_{n_1 | c} \sum_{n_2 = 1}^\infty \frac {A_F(n_2, n_1)} {n_1 n_2} S\left(a, - n_2; c / {n_1} \right) \Psi_-\left( \frac {n_2 n_1^2} {c^3} \right).
\end{split}
\end{equation*}
\end{prop}

\begin{rem}
We have used the functional equation 
\begin{equation*} 
\Gamma (1 - s) \Gamma (s) = \frac {\pi} {\sin (\pi s)}
\end{equation*}
and the fact that $k$ is even to rewrite the Gamma factors in \cite[Theorem 1.18]{MS} into the form in {\rm (\ref{eq: definition G eta (s)})}.
\end{rem}

\subsection{Gamma function}
Fix $s_0 \in \BC$ and let $ \Re \ s > - \Re \ s_0$. We have an asymptotic expansion as $|\Im \ s| \ra \infty$ 
\begin{equation*} 
\log \Gamma (s_0 + s) = \left( s_0 + s - \frac 1 2 \right) \log s - s + \frac 1 2 \log (2 \pi) + O\left( \frac 1 {|s| }\right).
\end{equation*}
With some calculations we show
\begin{equation}\label{eq: Stirling bound for G eta (s)}
\left| G_{\eta} (\sigma + i t) \right| \lll \left( {( |t| + 1) (t^2 + k^2)} \right)^{\sigma + \frac 1 2}.
\end{equation}

\subsection{Integral transform} \label{sec: integral transform}

Let $X > 0$. Suppose function $w(y)$ satisfies
\begin{equation}
\label{eq: conditions on the weight function}
\left\{ 
\begin{split}
& w(y) \text { is smooth with support in the dyadic interval } [X, 2X],\\
& y^j |w^{(j)}(y) | \leqslant c_j,
\end{split}\right.
\end{equation} for all $j \geqslant  0$ and some positive real numbers $c_j $. We call $w(y)$ an \textit{$X$-dyadic weight function}.

We are interested in $\psi(y ) = J_{2\kappa - 1} (\theta \sqrt { y}) w (y)$ with $\theta > 0$, $\kappa \geqslant 1$
and $w(y)$ an $X$-dyadic weight function with $y^j w^{(j)} (y) \lll_j 1$.

\subsubsection{ } \label{sec: theta sqrt X < 1} For $\theta \sqrt X \lll 1$ we use the bound (\ref{eq: bound for Bessel function}) for $J_{2\kappa - 1}$, then it follows from repeating integration by parts that
\begin{equation*}
\widetilde {\psi} (- \sigma - i t) \lll_{\sigma, A, \kappa} ( \theta \sqrt X)^{2\kappa - 1} X^{- \sigma} (|t| + 1)^{- A}.
\end{equation*}
This bound with $A = j + 3 \sigma + 3$ and the bound (\ref{eq: Stirling bound for G eta (s)}) for $G_\eta (s)$ yield, for $\sigma \geqslant  - \frac 1 2$,
\begin{equation*}
\begin{split}
 & \quad \quad x^j \Psi^{(j)}_{\pm} (x) \\
& \lll_{j, \sigma, \kappa} ( \theta \sqrt X)^{2\kappa - 1} \int_{-\infty}^\infty (x X)^{- \sigma} \left( (|t| + 1) (t^2 + k^2) \right)^{\sigma + \frac 1 2} \left( |t| + 1 \right)^{ - 3 \sigma - 3} dt\\
& \lll ( \theta \sqrt X)^{2\kappa - 1} k \left(\frac { k^2} {x X}\right)^{\sigma}.
\end{split}
\end{equation*}
If $x X \leqslant k^2 (k X)^\varepsilon$, we choose $\sigma = 0$ and obtain
\begin{equation}\label{eq: bound of xj Psi(j) pm}
x^j \Psi^{(j)}_{\pm} (x)
\lll_{j} ( \theta \sqrt X)^{2\kappa - 1} k.
\end{equation} Otherwise it is negligible by choosing $\sigma$ sufficiently large.

\subsubsection{ } \label{sec: theta sqrt X > 1} For $\theta \sqrt X \ggg 1$ we use the expression (\ref{eq: Bessel function and Whittaker function}) to write 
\begin{equation*}
\widetilde {\psi} (- \sigma - i t) = \sum_{\pm}  \frac 1 { \sqrt {2 \pi \theta} }  \int_0^\infty e^{\pm i \theta \sqrt y} W_{2\kappa - 1, \pm} (\theta \sqrt y) w (y) y^{-\sigma - i t - \frac 1 4} \frac {d y} y,
\end{equation*}
then repeating integration by parts and  the bound (\ref{eq: bound of Whittaker functions}) for $W_{ \pm, 2\kappa - 1}$ yield
\begin{equation*}
\widetilde {\psi} (- \sigma - i t) \lll_{\sigma, A, \kappa} \frac {1} {\sqrt \theta X^{\frac 1 4}} X^{-\sigma} \left(\frac {|t|} {\theta \sqrt X} + 1 \right)^{- A}.
\end{equation*}
This bound with $A = 3 \sigma + 3$ and the bound (\ref{eq: Stirling bound for G eta (s)}) for $G_\eta (s)$ yield, for $\sigma \geqslant  - \frac 1 2$,
\begin{equation*}
\begin{split}
& \quad \quad \Psi_{\pm} (x) \\
& \lll_{j, \sigma} \frac {1} {\sqrt \theta X^{\frac 1 4}} \int_{-\infty}^\infty (x X)^{- \sigma} \left( (|t| + 1) (t^2 + k^2) \right)^{\sigma + \frac 1 2} \left(\frac {|t|} { \theta \sqrt X} + 1 \right)^{ - 3 \sigma - 3} dt\\
& \lll \theta \sqrt X k \left(\frac {( \theta \sqrt X)^{3 } k^2} {x X}\right)^{\sigma}.
\end{split}
\end{equation*}
If $x X \leqslant ( \theta \sqrt X)^{3 } k^2 (k X)^\varepsilon$, we choose $\sigma = 0$ and obtain
\begin{equation}\label{eq: bound of Psi pm}
\Psi_{\pm} (x)
\lll_{j} \theta \sqrt X k.
\end{equation} Otherwise it is negligible by choosing $\sigma$ sufficiently large. 


\subsection{A Wilton-type bound}
We have the following Wilton-type bound involving conductor for $\Sym^2 f$ (\cite[Theorem 4.1]{Godber}). This type of bound was first proved for symmetric square lifts of $\GL(2, \BZ)$-Maass forms in \cite{LY}.

\begin{prop}\label{prop: Wilton-type bound}
Let $X > 0$ and $w(x)$ be an $X$-dyadic weight function defined in {\rm (\ref{eq: conditions on the weight function})}. Then for any real number $\alpha$,
\begin{equation}\label{eq: Wilton-type bound}
\sum_{n = 1}^\infty A_F(1, n) e(\alpha n) w(n) \lll_{\varepsilon, \{ c_j\}_j} X^{\frac 3 4 + \varepsilon} k^{\frac 1 2 + \varepsilon}.
\end{equation}
\end{prop}

\section{Proof of Theorem \ref{mthm}}
\label{sec: proof of Theorem}

\subsection{Amplified first moment average} 


Let $f \in S_k (1)$, $g_0 \in H^\star_{2\kappa} (P)$, and assume that $\kappa \geqslant  2$ is fixed, $ k > \kappa$ and $P$ is prime. Then $Q := k^4 P^3$ is essentially the conductor of $\Sym^2 f \otimes g$. 
We shall estimate the twisted first moment,
\begin{equation*} 
\mathfrak{F}_F(\ell ) := \sum_{g\in H^\star_{2\kappa} (P)}\omega_g^{-1}\lambda_g(\ell   ) L\left(\tfrac{1}{2},\Sym^2 f \otimes g \right)
\end{equation*}
when $(\ell,P)=1$ and $\ell \leqslant 16 L^4$ with $L$ to be chosen.  Subsequently, in the interest of simplifying notation, we shall write $F$ for $\Sym^2 f$.

Applying the approximate functional equation (\ref{eq: Approximate functional equation}), we have for $Y > 0$
\begin{equation*}
\begin{split}
L  &\left(\tfrac 1 2, F \otimes g \right) = \sum_{(d, P) = 1} \sum_{n = 1}^\infty \frac { \mu(d) A_F(1, n)  \lambda_g(d n)\lambda_g(\ell   )} { \sqrt {d^{ 3 } n }} V\left( \frac {d^3 n} {Y P^{\frac 3 2}}\right) \\
& \hskip 30 pt + (-1)^\kappa \sqrt P \sum_{(d, P) = 1} \sum_{n = 1}^\infty \frac { \mu(d) A_F(1, n) \lambda_g(d n P) \lambda_g(\ell   )} { \sqrt {d^{ 3 } n }} V\left( \frac {d^3 n Y} {P^{\frac 3 2}}\right),
\end{split}
\end{equation*}
where we have the applied multiplicative relation (\ref{eq: Hecke eigenvalues multiplicative 2}) which yields $\lambda_g(d n) \lambda_g( P) = \lambda_g(d n P)$.
Therefore, 
\begin{equation*}
\mathfrak F _F (\ell   ) = S_1 ( \ell, Y) + (-1)^\kappa S_2( \ell, Y),
\end{equation*}
with
\begin{equation*} 
\begin{split}
S_1( \ell   ) = S_1 ( \ell, Y) :=  \sum_{(d, P) = 1} \frac {\mu(d)} { d^{\frac  3 2}} \sum_{n =1}^\infty \frac { A_F(1, n)\Delta^\star_{2\kappa, P} (\ell   , dn)} { \sqrt n} V\left( \frac {d^3 n} {Y P^{\frac 3 2}}\right)
\end{split}
\end{equation*}
and
\begin{equation*} 
\begin{split}
S_2( \ell) = S_2( \ell, Y) := \sqrt P  \sum_{(d, P) = 1} \frac {\mu(d)} { d^{\frac  3 2}} \sum_{n =1}^\infty \frac { A_F(1, n) \Delta^\star_{2\kappa, P} (\ell   , dn P)} { \sqrt n } V\left( \frac {d^3 n Y} {P^{\frac 3 2}} \right).
\end{split}
\end{equation*}
It follows from the bound (\ref{eq: V(y)}) for $V$ that the contribution from $d^3 n > Y Q^{\frac 1 2 + \varepsilon}$ to $S_1( \ell   )$ and that from $d^3 n > \frac { Q^{\frac 1 2 + \varepsilon} } Y$ to $S_2( \ell   )$ are negligible.

Furthermore, we shall make the following \textit{a priori} assumption on $L$ and $Y$,
\begin{equation}\label{eq: a priori assumption}
{ L^2 \sqrt {Y} Q^{\frac 1 4} }  \leqslant P
\end{equation}
and see that our final choices of $Y$ in \eqref{Ychoice} and $L$ in \eqref{length-amp} satisfy this assumption.

\begin{rem}\label{rem: a priori assumption}
The assumption \eqref{eq: a priori assumption} will not be used until the Wilton-type bound (Proposition \ref{prop: Wilton-type bound}) is applied to the final estimates in \S \ref{sec: final estimate 1}. One reason of making this assumption is so that the weight function after the integral transform in Vorono\"i (see \S \ref{sec: theta sqrt X < 1}) satisfies \eqref{eq: conditions on the weight function} the hypothesis in Proposition \ref{prop: Wilton-type bound}.
\end{rem}




\subsection{Preparations for application of Petersson's trace formula}

\subsubsection{Treating $S_1( \ell)$ } Applying (\ref{eq: Petersson for Delta star}), 
the sum $S_1( \ell)$ is converted to
\begin{equation*}
\begin{split}
 \sum_{(d, P) = 1}  & {\mu(d)} \Bigg(\sum_{n = 1}^\infty \frac { A_F(1, n)} { \sqrt {d^3 n}} V\left( \frac {d^3 n} {Y P^{\frac 3 2}}\right) \Delta_{2\kappa, P} (\ell   , dn) \\
- & \sum_{n = 1}^\infty \frac { A_F(1, n)} { \sqrt {d^3 n}} V\left( \frac {d^3 n} {Y P^{\frac 3 2}}\right) \cdot \frac 1 {P v((n, P))} \sum_{j = 0}^{\infty} P^{- j} \Delta_{2\kappa, P} \lp \ell    P^{2j}, d n \rp \\
+ & \sum_{ n \equiv 0 (\mod P^2)} \frac { A_F(1, n)} { \sqrt {d^3 n}} V\left( \frac {d^3 n} {Y P^{\frac 3 2}}\right) \cdot \frac 1 {P+1} \sum_{j = 0}^{\infty} P^{- j} \Delta_{2\kappa, P} \lp \ell    P^{2j}, d n P^{-2} \rp \Bigg).
\end{split}
\end{equation*}
By trivial estimates using Deligne's bound, the bound (\ref{eq: V(y)}) for $V$ and the bound (\ref{eq: trivial bound for DeltakN}) for $\Delta_{2\kappa, P}$, the sum of the last two terms is bounded by
\begin{equation*}
\begin{split}
& \sum_{(d, P) = 1} \lp \sum_{ d^3 n \leqslant  Y Q^{\frac 1 2 + \varepsilon}} \frac {(\ell    d n P)^\varepsilon}  { \sqrt {d^3 n} P} +  \sum_{\sstyle n \equiv 0 (\mod P^2) \atop \sstyle d^3 n \leqslant  Y Q^{\frac 1 2 + \varepsilon} } \frac {(\ell    d n P)^\varepsilon}  { \sqrt {d^3 n} (P+1)} \rp \\
\lll & \frac {\sqrt Y Q^{\frac 1 4 + \varepsilon} } { P } \sum_{(d, P) = 1} \frac {1} { d^{ 3 - \varepsilon}} \lll \frac {\sqrt Y Q^{\frac 1 4 + \varepsilon}} { P }.
\end{split}
\end{equation*}
Therefore,
\begin{equation*}
S_1 ( \ell   ) = T_1 ( \ell   ) + O \lp \frac {\sqrt Y Q^{\frac 1 4  + \varepsilon} } { P } \rp,
\end{equation*}
with
\begin{equation} \label{eq: T1 (ell*)}
T_1 ( \ell   ) := \sum_{(d, P) = 1} {\mu(d)} \sum_{n = 1}^\infty \frac { A_F(1, n)} { \sqrt {d^3 n}}  V\left( \frac {d^3 n} {Y P^{\frac 3 2}}\right)  \Delta_{2\kappa, P} (\ell   , dn).
\end{equation}


\subsubsection{Treating $S_2 (\ell  )$} Similarly, after applying  (\ref{eq: Petersson for Delta star}), the sum $S_2 (\ell  )$ turns into
\begin{equation*}
\begin{split}
\sum_{(d, P) = 1}  & {\mu(d)}  \Bigg( \sqrt P \sum_{n = 1}^\infty \frac { A_F(1, n)} { \sqrt {d^3 n}} V\left( \frac {d^3 nY} { P^{\frac 3 2}}\right) \Delta_{2\kappa, P} \lp \ell   , dnP \rp \\
- & \sqrt P \sum_{n = 1}^\infty \frac { A_F(1, n)} { \sqrt {d^3 n}} V\left( \frac {d^3 nY} { P^{\frac 3 2}}\right) \cdot \frac 1 {P (P+1)}\sum_{j = 0}^{\infty} P^{- j} \Delta_{2\kappa, P} \lp P^{2j} \ell   , dnP \rp\\
+ & \sqrt P \sum_{n \equiv 0 (\mod P)} \frac { A_F(1, n)} { \sqrt {d^3 n}} V\left( \frac {d^3 nY} { P^{\frac 3 2}}\right) \cdot \frac 1 { P+1 }\sum_{j = 0}^{\infty} P^{- j} \Delta_{2\kappa, P} \lp P^{2j} \ell   , dnP\- \rp \Bigg).
\end{split}
\end{equation*}
By trivial estimates we have
\begin{equation*}
S_2 ( \ell   ) =  T_2 ( \ell   ) + O  \left(\frac {Q^{\frac 1 4 + \varepsilon}} {\sqrt Y P^{\frac 3 2} } \right),
\end{equation*}
with
\begin{equation} \label{eq: T2 (ell*)}
T_2 ( \ell   ) := \sqrt P \sum_{(d, P) = 1} {\mu(d)} \sum_{n = 1}^\infty \frac { A_F(1, n)} { \sqrt {d^3 n}} V\left( \frac {d^3 n Y } {P^{\frac 3 2}}\right) \Delta_{2\kappa, P} \lp \ell   , dnP \rp.
\end{equation}

\subsection{Application of Petersson's formula}

\subsubsection{Treating $T_1 (\ell  )$ } \label{sec: application of Petersson's formula 1} Applying Petersson's formula (\ref{eq: Petersson's formula}), $T_1 (\ell  )$ defined in (\ref{eq: T1 (ell*)}) is equal to
\begin{equation*}
\begin{split}
\sum_{(d, P) = 1} & {\mu(d)} \sum_{n = 1}^\infty \frac { A_F(1, n)} { \sqrt {d^3 n}}  V\left( \frac {d^3 n} {Y P^{\frac 3 2}}\right) \\
\times & \left( \delta (\ell   , dn) + (-1)^\kappa 2 \pi \sum_{ \scriptstyle c > 0 \atop \scriptstyle c \equiv 0 (\mod P) } \frac { S(\ell   , dn; c)} c J_{2\kappa-1} \left( \frac {4 \pi \sqrt { \ell   dn }} c \right) \right).
\end{split}
\end{equation*}
Trivial estimates using Deligne's bound and the bound (\ref{eq: V(y) y < k2}) for $V$ yields the following estimate on the diagonal term
\begin{equation*}
\sum_{\sstyle (d, P) = 1 \atop \sstyle d | \ell } \frac { A_F\left(1, \ell/d \right) } { \sqrt{\ell d^2} } V\left( \frac { \ell d^2} {  YP^{\frac 3 2}}\right) \lll \frac {Q^\varepsilon} {\sqrt \ell}.
\end{equation*}

Consider now the size of the off-diagonal terms
$$
\sum_{(d, P) = 1} {\mu(d)} \sum_{n = 1}^\infty \frac {A_F(1, n)} {\sqrt {d^3 n}}\sum_{ \scriptstyle c > 0 \atop \scriptstyle c \equiv 0 (\mod P) } \frac{S(\ell   , dn;c )}{c}J_{2\kappa-1}\left(\frac{4\pi\sqrt{ \ell   dn }}{c }\right) V\left( \frac {d^3 n} {Y P^{\frac 3 2}}\right).
$$
We wish to apply the Vorono\"i summation formula on the sum over  $n$. To this end, we subdivide the sum over $n$ by a smooth dyadic partition of unity for $V$, open the Kloosterman sum, and obtain the following sum up to a negligible error term due to (\ref{eq: V(y)}),
\begin{equation*} 
T_1^o (\ell  ) := \sum_{X \leqslant  Y Q^{\frac 1 2 + \varepsilon} } \sum_{(d, P) = 1}  {\mu(d)} T^{o}_1 (X, d, \ell  ),
\end{equation*}
where $X $ is of the form $2^{\frac j 2}$ with $j \geqslant  - 1$,
\begin{equation*} 
\begin{split}
 \quad \ T^{o}_1 (X, d, \ell  )  := \sum_{ \scriptstyle c > 0 \atop \scriptstyle c \equiv 0 (\mod P) } \frac 1 c \ \sumx_{a (\mod c )} e\left( \frac {\ell   a} c \right) \sum_{n = 1}^\infty {A_F(1, n)} e \left( \frac {d n \overline a} c \right) \psi_1 (n; c, X, d, \ell  ),
\end{split}
\end{equation*}
in which
\begin{equation*}
\psi_1 (y; c, X, d, \ell  ) := J_{2\kappa - 1} \left( \frac {4 \pi \sqrt { \ell   d y}} c \right) w_1 (y; X, d),
\end{equation*}
\begin{equation*}
w_1 (y; X, d) :=  \frac 1 { \sqrt{d^3 y} } V\left( \frac {d^3 y} { Y P^{\frac 3 2}}\right) h\left( \frac {d^3 y} {X}\right),
\end{equation*}
with $h (y)$ some smooth function supported on $[1, 2]$ satisfying $h^{(j)} (y) \lll_j 1.$ In view of the range of $X$, the bound (\ref{eq: V(y) y < k2}) for $V$ implies that $w_1  (y; X, d)$ is an $X$-dyadic weight function with bounds 
\begin{equation}\label{eq: bounds for w1}
x^j w^{(j)}_1  (y; X, d) \lll_{j, \varepsilon, l} \frac { Q^\varepsilon } {\sqrt X}.
\end{equation}

Moreover, we pull out the greatest common divisor $(c,d)$ by writing
\begin{equation*} 
T^{o}_1 (X, d, \ell  ) = \sum_{d_1 d_2 = d} \sum_{\sstyle (c, d_2) = 1 \atop \sstyle c \equiv 0 (\mod P)} T^o_1(c, X, d_1, d_2,  \ell  ),
\end{equation*}
with
\begin{equation} \label{eq: To1(c, X, d1, d2, ell*)}
\begin{split}
T^o_1(c, X, d_1, d_2, \ell  ) = \frac 1 {c d_1 }\  & \sumx_{a (\mod c d_1 P)} e\left( \frac {\ell   a } {c d_1 }\right) \\
& \sum_{n = 1}^\infty {A_F(1, n)} e \left( \frac {d_2 n \overline a} c \right) \psi_1 (n; c d_1, X, d_1 d_2,  \ell  ),
\end{split}
\end{equation}

\subsubsection{Treating $T_2( \ell  )$ } The diagonal term coming from Petersson's formula applied to $T_2( \ell  )$ vanishes since $(\ell,P)=1$. The off-diagonal sum is of size
\begin{equation*} 
\sum_{(d, P) = 1} {\mu(d)} \sum_{n = 1}^\infty  \frac { A_F(1, n)} { \sqrt {d^3 n}}  V\left( \frac {d^3 n} {Y P^{\frac 3 2}}\right) 
 \sum_{ c = 1 }^\infty \frac {S( \ell  , d n P; c P)} {c \sqrt P} J_{2\kappa-1} \left( \frac {4 \pi \sqrt { \ell   d n}} {c \sqrt P} \right).
\end{equation*}
Since $S( \ell  , d n P; c P) = 0$ when $P | c$, we may impose the condition $(c, P) = 1$ to the sum over $c$.
Under this condition
\begin{equation*}
S(\ell  , dn P; cP) = S(\ell   \overline c, 0; P) S(\ell   \overline P, d n; c) = - S(\ell   \overline P, d n; c).
\end{equation*}
Inserting this above and reordering our sums we obtain
\begin{equation*} 
\begin{split}
\sum_{(d, P) = 1} {\mu(d)} \sum_{(c, P) = 1} \frac {1} {c \sqrt P} \sum_{n = 1}^\infty \frac { A_F(1, n)} { \sqrt {d^3 n}}  S(\ell   \overline P, dn; c) J_{2\kappa - 1} \left( \frac {4 \pi \sqrt {\ell   d n}} {c \sqrt P} \right) V\left( \frac {d^3 n} {Y P^{\frac 3 2}}\right).
\end{split}
\end{equation*}

With the same treatment as in \S \ref{sec: application of Petersson's formula 1}, one is left to study the sum
\begin{equation*} 
T_2^o (\ell  ) := \sum_{X \leqslant  { Q^{\frac 1 2 + \varepsilon} }/ Y}  \sum_{(d, P) = 1}  {\mu(d)} T^{o}_2 (X, d, \ell  )
\end{equation*}
with
\begin{equation*}
T^{o}_2 (X, d, \ell  ) := \sum_{d_1 d_2 = d} \sum_{(c, P d_2) = 1} T^o_2 (c, X, d_1, d_2, \ell  ),
\end{equation*}
and
\begin{equation*} 
\begin{split}
T^{o}_2 (c, X, d_1, d_2, \ell  ) := \frac {1} {c d_1 \sqrt P } \  & \sumx_{a (\mod c d_1)} e\left( \frac {\ell   \overline P a} { c d_1 }\right)\\
& \sum_{n = 1}^\infty {A_F(1, n)} e \left( \frac {d_2 n \overline a} c \right) \psi_2 (n; c d_1, X, d_1 d_2, \ell  ),
\end{split}
\end{equation*}
where
\begin{equation*}
\psi_2 (y; c, X, d, \ell  ) :=  J_{2\kappa - 1} \left( \frac {4 \pi \sqrt { \ell   d y }} {c \sqrt P} \right) w_2 (y; X, d)
\end{equation*}
and
\begin{equation*}
w_2 (y; X, d) :=  \frac 1  {\sqrt {d^3 y} } V\left( \frac {d^3 y } {Y P^{\frac 3 2}}\right) h\left( \frac {d^3 y} {X}\right).
\end{equation*}
Again, $w_2  (y; X, d)$ is an $X$-dyadic weight function with bounds 
\begin{equation}\label{eq: bounds for w2}
x^j w^{(j)}_2  (y; X, d) \lll_{j, \varepsilon, l} \frac { Q^\varepsilon} {\sqrt X}.
\end{equation}

\subsection{Application of Vorono\"i's formula}

\subsubsection{ Treating $T^o_1 (c, X, d_1, d_2, \ell  )$ }\label{sec: application of Voronoi's formula 1} Applying Vorono\"i's formula to the innermost sum in (\ref{eq: To1(c, X, d1, d2, ell*)}),  $T^o_1 (c, X, d_1, d_2, \ell  )$ is converted to the sum
$$ T_{1, +}^o (c, X, d_1, d_2, \ell  ) + T_{1, -}^o (c, X, d_1, d_2, \ell  ),$$
where
\begin{equation*}
\begin{split}
T_{1, \pm}^o (c, X, & d_1, d_2, \ell  ) := \frac 1 {d_1}\ \sumx_{a (\mod c d_1)} e\left( \frac {\ell  a} {c d_1}\right)
\\ 
& \sum_{n_1 c_1 = c}  \sum_{n_2 = 1}^\infty \frac {A_F(n_2, n_1)} {n_1 n_2} S (a \overline {d}_2, \pm n_2; c_1 ) \Psi_{1, \pm} \left( \frac {n_2 n_1^2} {c^3} ; c d_1, X, d_1 d_2, \ell   \right),
\end{split}
\end{equation*}
and $\Psi_{1, \pm} (x; c, X, d, \ell  )$ is the integral transform of $\psi_{1, \pm} (y; c, X, d, \ell  )$ defined as in \S \ref{sec: Voronoi}. Opening the Kloosterman sum and changing the order of summation in $T_{1, \pm}^o (c, X, d_1, d_2, \ell  )$ above, we arrive at
\begin{equation*}
\begin{split}
T_{1, \pm}^o (c, X, d_1, d_2, \ell  ) 
 = \frac 1 {d_1}  & \sum_{n_1 c_1 = c} \frac 1 { n_1} \ \sumx_{b (\mod c_1)} \  \sumx_{a (\mod c d_1)} e\left( \frac {a(b d_1 n_1 \overline {d}_2 + \ell  )} {c d_1}\right) \\
& \sum_{n_2 = 1}^\infty \frac {A_F(n_2, n_1)} { n_2} e\left( \pm \frac {n_2 \overline b} {c_1}\right) \Psi_{1, \pm} \left( \frac {n_2 n_1^2} {c^3} ; c d_1, X, d_1 d_2, \ell   \right).
\end{split}
\end{equation*}
By M\"obius inversion, 
\begin{equation*}
\sumx_{a (\mod c d_1)} e\left( \frac {a (b d_1 n_1 \overline {d}_2 + \ell  )}  {c d_1} \right) = \sum_{c_2 | c d_1 } \mu \left( \frac {c d_1} {c_2}\right) \sum _{a (\mod c_2)} e\left( \frac {a (b d_1 n_1 \overline {d} _2+ \ell  )} {c_2}\right).
\end{equation*}
The inner sum produces a congruence condition $b d_1 n_1 \equiv - d_2 \ell  (\mod c_2)$, which forces $(c_2, d_1 n_1 ) | (c_2, \ell  )$ and hence $c_2 | c_1 ( d_1 n_1, \ell  )$. Therefore, $T_{1, \pm}^o (X, d_1, d_2, \ell)$ becomes
\begin{equation*}
\begin{split}
\frac 1 {d_1} & \sum_{n_1 c_1 = c}\ \sum_{ \scriptstyle c_2 | c_1 (d_1 n_1, \ell  ) } \mu \left( \frac {c d_1} {c_2}\right) \frac { c_2 } {n_1} \\ 
& 
\sumx_{\sstyle b (\mod c_1) \atop \sstyle b \equiv - \ell d_2 \overline { d_1 n_1} (\mod c_2)} \sum_{n_2 = 1}^\infty e\left( \pm \frac {n_2 \overline b} {c_1}\right) \frac {A_F(n_2, n_1)} { n_2} \Psi_{1, \pm} \left( \frac {n_2 n_1^2} {c^3} ; c d_1, X, d_1 d_2, \ell   \right).
\end{split}
\end{equation*}
Finally, we apply the Hecke relation (\ref{eq: Hecke relation}),
\begin{equation} \label{eq: T1pmo (X; d_1, d_2, ell)}
\begin{split}
&\ T_{1, \pm}^o (c, X, d_1, d_2, \ell ) \\
= & \frac 1 {d_1} \sum_{n_1 c_1 = c} \ \sum_{ \scriptstyle c_2 | c_1 (d_1 n_1, \ell  ) } \sum_{n_3 | n_1} \mu \left( \frac {c d_1} {c_2}\right) \mu (n_3) \frac { c_2} {n_1 n_3}  A_F(1, n_1/n_3 )\\ 
& 
\sumx_{\sstyle b (\mod c_1) \atop \sstyle  b \equiv - \ell d_2 \overline { d_1 n_1}   (\mod c_2)} \sum_{n = 1}^\infty \frac {A_F(n, 1)} { n }  e\left( \pm \frac {n n_3 \overline b} {c_1}\right) \Psi_{1, \pm} \left( \frac {n n_3 n_1^2} {c^3} ; c d_1, X, d_1 d_2, \ell   \right).
\end{split}
\end{equation}

\subsubsection {Treating $T^o_2 (c, X, d_1, d_2, \ell  )$ }  Following the same line of arguments as in \S \ref{sec: application of Voronoi's formula 1}, we have
$$T^o_2 (c, X, d_1, d_2, \ell  ) = T_{2, +}^o (c, X, d_1, d_2, \ell  ) + T_{2, -}^o (c, X, d_1, d_2, \ell  ) $$
with
\begin{equation} \label{eq: T2pmo (X; d_1, d_2, ell)}
\begin{split}
& \ T_{2, \pm}^o (c, X; d_1, d_2, \ell  ) \\
= & \frac 1 {\sqrt P d_1} \sum_{n_1 c_1 = c} \ \sum_{ \scriptstyle c_2 | c_1 (d_1 n_1, \ell  ) } \sum_{n_3 | n_1} \mu \left( \frac {c d_1} {c_2}\right) \mu (n_3) \frac { c_2} {n_1 n_3}  A_F(1, n_1/n_3 )\\ 
& \sumx_{\sstyle b (\mod c_1) \atop \sstyle b \equiv - \ell   d_2 \overline { d_1 n_1 P} (\mod c_2)}
\sum_{n = 1}^\infty \frac {A_F(n, 1)} { n } e\left( \pm \frac {n n_3 \overline b} {c_1}\right) \Psi_{2, \pm} \left( \frac {n n_3 n_1^2} {c^3} ; cd_1, X, d_1 d_2, \ell   \right).
\end{split}
\end{equation}

\subsection{Final estimates}
We restrict ourselves to the partial sums of (\ref{eq: T1pmo (X; d_1, d_2, ell)}) and (\ref{eq: T2pmo (X; d_1, d_2, ell)}) in which $d = n_1 = 1$ to simplify the complicated notation; the general cases may be treated in the same way.

\subsubsection{ } \label{sec: final estimate 1}
Suppose $X \leqslant  Y { Q^{\frac 1 2 + \varepsilon} }$ and $c \geqslant P$ with $c \equiv 0 (\mod P)$. Let
\begin{equation*}
\psi_1 ( y; c, X, \ell   ) := \psi_1  ( y; c, X, 1, \ell   ) = J_{2\kappa - 1} \left( \frac {4 \pi \sqrt { \ell  y }} c \right) w_1 (y; X),
\end{equation*}
with 
\begin{equation*}
w_1 (y; X) := w_1 (y; X, 1) = \frac 1 {  \sqrt y} V\left( \frac { y Y} {P^{\frac 3 2}}\right)  h\left( \frac { y} {X}\right).
\end{equation*}
We consider the sum
\begin{equation} \label{eq: T1pmo (c, X, ell)}
\begin{split}
T_{1, \pm}^o (c, X, \ell  ) &
:= \sum_{\sstyle c_2 | c \atop \sstyle (c_2, \ell) = 1}  \mu \left( \frac {c } {c_2}\right) { c_2}\\
& \sumx_{\sstyle b (\mod c ) \atop \sstyle b \equiv - \ell (\mod c_2)} \sum_{n = 1}^\infty \frac {A_F(n, 1)} { n } e\left( \pm \frac {n \overline b} {c }\right) \Psi_{1, \pm} \left( \frac {n } {c^3} ; c, X, \ell   \right).
\end{split}
\end{equation}
Then summing all $ T_{1, \pm}^o (c, X, \ell  )$ over $c$ and $\pm $ yields the partial sum of $T_1^o (X, 1, \ell)$ with $n_1 = 1$.

By our assumption (\ref{eq: a priori assumption}) we have $\frac {4 \pi \sqrt { \ell X}} c \lll \frac { L^2 \sqrt X} P \leqslant 1$, so we are in the situation of \S \ref{sec: theta sqrt X < 1} with $\theta = \frac {4 \pi \sqrt { \ell}} c$ and $w (y) = w_1 (y; X) \frac { \sqrt X} { Q^{ \varepsilon} }$ in view of (\ref{eq: bounds for w1}). According to (\ref{eq: bound of xj Psi(j) pm}), for $x X \leqslant  k^2 Q^\varepsilon$ we have the following bounds
\begin{equation}\label{eq: bound of Psi 1 pm ( x ; c, X, ell)}
x^j \Psi_{1, \pm}^{(j)} \left( x ; c, X, \ell   \right) \lll_{j, \varepsilon} \frac {k Q^\varepsilon} {\sqrt X } \left(\frac {\sqrt {\ell X }} c \right)^{2\kappa - 1}.
\end{equation}
Otherwise, the bounds are arbitrarily small.

Next we wish to use the Wilton-type bound in the form (\ref{eq: Wilton-type bound}). For this we apply a dyadic partition of unity and convert the innermost sum over $n$ in (\ref{eq: T1pmo (c, X, ell)}) into
\begin{equation} \label{eq: sum Z additive twist}
\sum_{Z} \sum_{n = 1}^\infty {A_F(n, 1)} e\left( \pm \frac {n \overline b} {c }\right) w_{1, \pm} (n; c, X, Z,  \ell  ),
\end{equation}
where $Z = 2^{\frac j 2}$ with $j \geqslant  - 1$ and
\begin{equation*}
w_{1, \pm} (x; c, X, Z,  \ell  ) := \frac 1 x \Psi_{1, \pm} \left( \frac {x} {c^3} ; c, X, \ell   \right) h\left( \frac x {Z} \right).
\end{equation*} 
The contribution from $Z > \frac {c^3 k^2 Q^\varepsilon} X $ is negligible. Otherwise, the bound (\ref{eq: bound of Psi 1 pm ( x ; c, X, ell)}) implies that $ w_{1, \pm} (x; c, X, Z,  \ell)$ is a $Z$-dyadic weight function with 
$$x^j w_{1, \pm}^{(j)} (x; c, X, Z, \ell) \lll_{j, \varepsilon} \frac {k Q^\varepsilon} {\sqrt X Z} \left(\frac {\sqrt {\ell X}} c \right)^{2\kappa - 1}.$$
Applying the Wilton-type bound (\ref{eq: Wilton-type bound}) to the inner sum in (\ref{eq: sum Z additive twist}) gives
\begin{equation*}
\begin{split}
T_{1, \pm}^o (c, X, \ell ) & \lll \sum_{\sstyle c_2 | c \atop \sstyle (c_2, \ell) = 1} c \sum_{Z \leqslant  {c^3 k^2 Q^\varepsilon}/ X} \frac { k^{\frac 3 2} Q^{\varepsilon}} { \sqrt X Z^{\frac 1 4} } \left(\frac {\sqrt {\ell X}} c \right)^{2\kappa - 1} \\
& \lll \frac { c^{\frac 1 4 + \varepsilon} k Q^{\varepsilon}} {X^{\frac 1 4} }  \left(\frac {\sqrt {\ell X }} c \right)^{2\kappa - 1}.
\end{split}
\end{equation*}
Therefore,
\begin{equation}\label{eq: c-sum requiring kappa > 1, 1}
\sum_{\sstyle c \geqslant  P \atop \sstyle c \equiv 0 (\mod P)}  T_{1, \pm}^o (c, X, \ell ) \lll \frac { k Q^{\varepsilon}} {X^{\frac 1 4}} \sum_{\sstyle c \geqslant  P \atop \sstyle c \equiv 0 (\mod P)}  c^{\frac 1 4 + \varepsilon} \left(\frac {\sqrt {\ell X}} c \right)^{2\kappa - 1} \lll \frac {\ell^{\frac 5 8} X^{\frac 3 8} k Q^{\varepsilon}} P,
\end{equation}
where it should be noted that the assumption $\kappa \geqslant 2$ is required to guarantee the convergence of the $c$-sum.

Finally, summing over $X$, the partial sum of $T_1^o (\ell)$ in which $d = n_1 = 1$ is bounded by
\begin{equation*}
\sum_{X \leqslant  Y Q^{\frac 1 2 + \varepsilon}} \frac {\ell^{\frac 5 8} X^{\frac 3 8}  k Q^{\varepsilon}} P \lll \frac {\ell^{\frac 5 8} Y^{\frac 3 8}  k Q^{\frac 3 { 16} + \varepsilon}} P.
\end{equation*}


\subsubsection{ } 
Suppose $X \leqslant   \frac { Q^{\frac 1 2 + \varepsilon} } Y$ and $(c, P) = 1$. Let
\begin{equation*}
\psi_2 ( y ; c, X, \ell ) := \psi_2 ( y; c, X, 1, \ell ) = J_{2\kappa - 1} \left( \frac {4 \pi \sqrt {\ell y }} {c \sqrt P }\right) w_2 (y; X),
\end{equation*}
with 
\begin{equation*}
w_2 (y; X) := w_2 (y; X, 1) = \frac 1 {\sqrt y} V\left( \frac {y} { Y P^{\frac 3 2}}\right) h\left( \frac {y} {X}\right).
\end{equation*}
Consider the sum
\begin{equation}  \label{eq: T2pmo (c, X, ell)}
\begin{split}
T_{2, \pm}^o (c, X, \ell) &  := \frac 1 {\sqrt P} \sum_{\sstyle c_2 | c \atop \sstyle (c_2, \ell) = 1}  \mu \left( \frac {c } {c_2}\right)  { c_2}\\
& \sumx_{\sstyle b (\mod c ) \atop \sstyle b \equiv - \ell \overline P (\mod c_2)}
\sum_{n = 1}^\infty \frac {A_F(n, 1)} { n } e\left( \pm \frac {n \overline b} {c }\right) \Psi_{2, \pm} \left( \frac {n } {c^3} ; c , X, \ell \right).
\end{split}
\end{equation}
All $ T_{2, \pm}^o (c, X, \ell  )$ constitute the partial sum of $T_2^o (X, 1, \ell)$ with $n_1 = 1$.

We shall apply arguments similar to those in \S \ref{sec: final estimate 1} except that, instead of using the Wilton-type bound, Deligne's bound is trivially applied to each individual $A_F(n, 1)$. 

When $c \geqslant  \sqrt {\frac {\ell X} P}$, the analysis of \S \ref{sec: theta sqrt X < 1} implies that we may truncate the sum over $n$ in (\ref{eq: T2pmo (c, X, ell)}) at $\frac {c^3 k^2 Q^\varepsilon} X$ with the difference of a negligible error, and for $ n \leqslant  \frac {c^3 k^2 Q^\varepsilon} X$ we have
\begin{equation*}
\Psi_{2, \pm} \left( \frac {n } {c^3} ; c , X, \ell \right) \lll \frac { k Q^{\varepsilon}} { \sqrt X } \left(\frac {\sqrt {\ell X}} {c \sqrt P} \right)^{2\kappa - 1}.
\end{equation*}
Thus trivial estimates yield
\begin{equation}\label{eq: T 2pmo (c, X, ell) 1}
\begin{split}
T_{2, \pm}^o (c, X, \ell) & \lll \frac 1 {\sqrt P} \sum_{ c_2 | c } c \sum_{n \leqslant  {c^3 k^2 Q^\varepsilon}/ X} \frac { k Q^{\varepsilon}} { \sqrt X n} \left(\frac {\sqrt {\ell X}} {c \sqrt P} \right)^{2\kappa - 1} \\
& \lll \frac { c^{1 + \varepsilon} k Q^{\varepsilon}} {\sqrt { X P}} \left(\frac {\sqrt  {\ell X}} {c \sqrt P} \right)^{2\kappa - 1}.
\end{split}
\end{equation}
When $c < \sqrt {\frac {\ell X} P}$, we are in the situation of \S \ref{sec: theta sqrt X > 1}. It follows from \eqref{eq: bound of Psi pm} that if $n$ does not exceed $\frac {\left(\sqrt {\ell X}/(c\sqrt P) \right)^3 c^3  k^2 Q^\varepsilon} {X}$ then we have the bound
\begin{equation*}
\Psi_{2, \pm} \left( \frac {n } {c^3} ; c , X, \ell \right) \lll  \frac { k Q^{\varepsilon}} { \sqrt X } \frac {\sqrt  {\ell X}} {c \sqrt P}.
\end{equation*}
Otherwise, we have a negligible contribution.
Therefore,
\begin{equation}\label{eq: T 2pmo (c, X, ell) 2}
\begin{split}
T_{2, \pm}^o (c, X, \ell) & \lll \frac 1 {\sqrt P} \sum_{ c_2 | c } c \sum_{n \leqslant  {\left(\sqrt  {\ell X}/(c\sqrt P) \right)^3 c^3  k^2 Q^\varepsilon} \textstyle {/} \sstyle  X} \frac { k Q^{\varepsilon}} { \sqrt X n} \frac {\sqrt  {\ell X}} {c \sqrt P} \\
& \lll \frac { c^{1 + \varepsilon} k Q^{\varepsilon}} {\sqrt {X P }} \frac {\sqrt {\ell X}} {c \sqrt P}.
\end{split}
\end{equation}

Combining (\ref{eq: T 2pmo (c, X, ell) 1}, \ref{eq: T 2pmo (c, X, ell) 2}) we have
\begin{equation} \label{eq: c-sum requiring kappa > 1, 2}
\begin{split}
& \sum_{ (c, P) = 1 } T_{1, \pm}^o (c, X, \ell )\\
\lll &  \frac { k Q^{\varepsilon}} {\sqrt { X P}} \left( \sum_{\sstyle c \geqslant  \sqrt {\ell X/P} \atop \sstyle (c, P) = 1 }  c^{1 + \varepsilon} \left(\frac {\sqrt {\ell X}} {c \sqrt P} \right)^{2\kappa - 1} +  \sum_{\sstyle c < \sqrt {\ell X/P} \atop \sstyle (c, P) = 1 }   c^{1 + \varepsilon} \frac {\sqrt {\ell X}} {c \sqrt P} \right) \lll \frac {\ell \sqrt { X}  k Q^{\varepsilon} } {P^{\frac 3 2}}.
\end{split}
\end{equation}
Again, the assumption $\kappa \geqslant 2$ guarantees the convergence of this $c$-sum.

We conclude with the following bound for the partial sum of $T_2^o (\ell)$ where $d = n_1 = 1$,
\begin{equation*} 
\begin{split}
\sum_{X \leqslant  Q^{\frac 1 2 + \varepsilon}/Y} \frac {\ell \sqrt { X}  k Q^{\varepsilon} } {P^{\frac 3 2}} \lll \frac {\ell k  Q^{\frac 1 4 + \varepsilon}} {\sqrt Y P^{\frac 3 2} }.
\end{split}
\end{equation*}


\subsection{Conclusion}
In conclusion, summing all contributions in the above arguments, we have the following bound for the twisted first moment
\begin{equation*}
\begin{split}
\mathfrak{F}_F(\ell ) & = \sum_{g\in H^\star_{2\kappa} (P)} \omega_g^{-1} \lambda_g(\ell) L \left(\tfrac{1}{2}, F \otimes g\right)\\
& \lll \left( \frac{\sqrt{Y} Q^{\frac{1}{4}}}{P} + \frac {Q^{\frac 1 4 }} {\sqrt Y P^{\frac 3 2} } + \frac{1}{\sqrt{\ell}} + \frac {\ell^{\frac{5}{8}} Y^{\frac{3}{8}}kQ^{\frac{3}{16}}}{P} + \frac {\ell k Q^{\frac{1}{4}}} { \sqrt{Y} P^{\frac 3 2}} \right) 
Q^{\varepsilon},
\end{split}
\end{equation*}
where $\ell \leqslant 16 L^4$ is co-prime with $P$ and with $L$ and $Y$ satisfying our assumption \eqref{eq: a priori assumption}.
Thus,
\begin{equation}
\label{eq: bd-for-twisted-mom}
\begin{split}
\mathfrak{F}_F(\ell ) \lll & \left(\frac 1 {\sqrt \ell} + \frac {\sqrt Y Q^{\frac 1 4 } } { P } + \frac { Q^{\frac 1 4 }} {\sqrt Y P^{\frac 3 2} } + \frac {L^{\frac 5 2} Y^{\frac 3 8}  k Q^{\frac 3 {16} }} P + \frac {L^4 k  Q^{\frac 1 4 }} { \sqrt Y P^{\frac 3 2}} \right) Q^{\varepsilon}\\
\lll & \left(\frac 1 {\sqrt \ell}  + \frac { \sqrt Y Q^{\frac 1 4 } } { P } +  \frac {L^{\frac 5 2} Y^{\frac 3 8}  k Q^{\frac 3 {16} }} P + \frac {L^4 k  Q^{\frac 1 4 }} { \sqrt Y P^{\frac 3 2}} \right) Q^{\varepsilon}\\
\lll & \left(\frac 1 {\sqrt \ell}  + \frac {L^{\frac 5 2} Y^{\frac 3 8}  k Q^{\frac 3 {16} }} P + \frac {L^4 k  Q^{\frac 1 4 }} { \sqrt Y P^{\frac 3 2}} \right) Q^\varepsilon,
\end{split}
\end{equation}
where the last line follows from our assumption $ { L^2 \sqrt {Y} Q^{\frac 1 4} } \leqslant P$ in \eqref{eq: a priori assumption}.
We achieve an optimal bound by choosing 
\begin{equation}\label{Ychoice}
Y = L^{\frac {12}7} k^{\frac 2 7} P^{-\frac 5 {14}}. 
\end{equation}
Inserting this value for $Y$ in \eqref{eq: bd-for-twisted-mom} we get Theorem~\ref{mthm}. 




\bibliographystyle{plain}
\bibliography{references}

\end{document}